\magnification=\magstep1
\hsize16truecm
\vsize23.5truecm
\topskip=1truecm
\raggedbottom
\abovedisplayskip=3mm
\belowdisplayskip=3mm
\abovedisplayshortskip=0mm
\belowdisplayshortskip=2mm
\normalbaselineskip=12pt
\normalbaselines
\font\titlefont= cmcsc10 at 12pt
\def\F{\Bbb F}

\def\C{\Bbb C}
\def\Z{\Bbb Z}
\def\Q{\Bbb Q}
\def\P{\Bbb P}
\def\A{\Bbb A}
\def\G{\Bbb G}

\def\E{\Bbb E}

\def\tto{\longrightarrow}

%
%
\catcode`\@=11
\font\tenmsa=msam10
\font\sevenmsa=msam7
\font\fivemsa=msam5
\font\tenmsb=msbm10
\font\sevenmsb=msbm7
\font\fivemsb=msbm5
\newfam\msafam
\newfam\msbfam
\textfont\msafam=\tenmsa  \scriptfont\msafam=\sevenmsa
  \scriptscriptfont\msafam=\fivemsa
\textfont\msbfam=\tenmsb  \scriptfont\msbfam=\sevenmsb
  \scriptscriptfont\msbfam=\fivemsb
\def\hexnumber@#1{\ifcase#1 0\or1\or2\or3\or4\or5\or6\or7\or8\or9\or
	A\or B\or C\or D\or E\or F\fi }
\edef\msa@{\hexnumber@\msafam}
\edef\msb@{\hexnumber@\msbfam}
\mathchardef\square="0\msa@03
\mathchardef\subsetneq="3\msb@28
\mathchardef\ltimes="2\msb@6E
\mathchardef\rtimes="2\msb@6F
\def\Bbb{\ifmmode\let\next\Bbb@\else
	\def\next{\errmessage{Use \string\Bbb\space only in math mode}}\fi\next}
\def\Bbb@#1{{\Bbb@@{#1}}}
\def\Bbb@@#1{\fam\msbfam#1}
\catcode`\@=12
%
%
%
\noindent
\font\eighteenbf=cmbx10 scaled\magstep3
\font\titlefont=cmcsc10 at 12pt
\vskip 5.0pc
\centerline{\eighteenbf  The Chow Ring of the  Moduli
Space}\footnote{}{This version corrects two inaccuracies of the 1998 JAG paper}
\bigskip
\centerline{\eighteenbf  of Abelian Threefolds}
\vskip 2pc
\centerline{\titlefont Gerard van der Geer}
\vskip2pc
\noindent In this paper we  determine the structure of the Chow ring of the
Delaunay-Voronoi compactification $\tilde{\cal A}_3$ of the moduli space of
principally polarized abelian threefolds. This compactification was
introduced by Namikawa and studied by Tsushima. We shall use equivariant
classes on  level coverings of $\tilde{\cal A}_3$. We also compare this ring
with the Chow ring of the moduli space of stable genus $3$ curves as
determined by Faber. 
\vskip 2.0pc
\noindent
{\bf \S 0 Introduction}
\bigskip
\noindent
Let ${\cal A}_g$ be the moduli stack of principally polarized abelian
varieties of dimension $g$ and  ${\cal X}_g$ the universal family of abelian
varieties over ${\cal A}_g$. The stack ${\cal A}_g$ is not complete; it can
be completed to a space ${\cal A}_g^*$, the so-called Satake
compactification. In Faltings-Chai [F-C] smooth compactifications of ${\cal
A}_g$ are constructed. These compactifications are not unique as they depend
on the choice of a polyhedral decomposition of the cone of positive definite
symmetric bilinear forms in $g$ real variables. However, a specific choice
of polyhedral decomposition (the Delaunay-Voronoi decomposition) gives for
low values of $g$ a  good (smooth) compactification, cf.\ [N].
For $g=1$ it gives the  unique compactification, while for  $g=2$ one 
gets the blow-up of the Satake compactification as considered by Igusa.  For
$g=3$ one gets the  compactification $\tilde{\cal A}_3$ as constructed by
Namikawa and studied by Tsushima in [N] and [T 1,2]. This is a smooth stack
and it has a Chow ring. In the following $\tilde{\cal A}_g$ will denote this
Delaunay-Voronoi compactification.
\par
The Chow ring of $\tilde{\cal A}_2$ was determined by Mumford in [M3],
since he determined the Chow ring of $\overline{\cal M}_2$. In this paper we
determine the structure of the Chow ring of $\tilde{\cal A}_3$. It is a ring
with four generators.  We do this by bounding the ranks of the Chow
groups and then use the intersection numbers of equivariant classes to get
the generators. In [F] Faber determined the Chow ring of ${\overline {\cal
M}}_3$. For our choice of compactification the Torelli morphism extends
and we can use this  morphism $t: {\overline {\cal M}}_3 \to \tilde{\cal
A}_3$ of degree $2$ to compare both rings. 
\par
I like to thank Carel Faber and Valery Alexeev for some helpful remarks.
I also thank Sam Grushevsky for pointing out an inaccuracy.
\vskip 1.0pc
\noindent
{\bf \S 1 Preliminaries}
\bigskip
\noindent
The Satake compactification ${\cal A}_g^*$ of the stack ${\cal A}_g$
over $\Z$ admits a stratification
$$
{\cal A}_g^*= \coprod_{j=0}^g {\cal A}_j.
$$
Let ${\cal A}_g[\ell]$ denote the moduli stack over $\Z[1/\ell]$ of
principally polarized abelian varieties of dimension $g$ with a full
symplectic level
$\ell$ structure. There is a morphism $p_{\ell}: {\cal A}_g[\ell] \tto {\cal
A}_g$ which is equivariant for the natural action of $G_{\ell}:=Sp(2g,\Z /
\ell \Z)$ on
${\cal A}_g[\ell]$. The degree of this morphism is
$$
\gamma=\gamma_{\ell} = \ell^{g(2g+1)}\prod_{p|\ell} \prod_{1\leq j \leq
g}(1-p^{-2j}).
$$
The morphism $p_\ell$ extends to the Satake compactification
${\cal A}_g^*[\ell]$. For $\ell \geq 3$ we denote by 
$$
\mu_g= \mu_g(\ell)= {1 \over 2} \ell^{2g} \prod_{p|\ell}\prod_{1\leq j \leq
g} (1-p^{-2j})
$$
the number of maximal dimensional cusps of ${\cal A}_g^*[\ell]$ (each of
them is a copy of ${\cal A}_{g-1}^*[\ell]$). 
\par
The toroidal compactifications we need are described by Namikawa in [N] and
by Faltings-Chai [F-C], see also [T1]. There one also finds the toroidal
compactification of the universal abelian variety, see [F-C, p. 195]. The
toroidal compactifications of Faltings-Chai  carry a rank
$g$ vector bundle, the Hodge bundle $\E$. It is given over ${\cal A}_g$ by
$s^*(\Omega_{X/S})$ for each principally polarized abelian variety
$X/S$ corresponding to $S\to {\cal A}_g$ with zero section $s$. The pull back
$\pi^*(\E)$ under $\pi: \tilde{\cal X}_g \tto \tilde {\cal A}_g$ of the Hodge
bundle to the universal scheme $\tilde{\cal X}_g$ over $\tilde{\cal A}_g$ is
isomorphic to the sheaf of relative logarithmic differentials
$\Omega_{\tilde{\cal X}}(d \log \infty)/\pi^*(\Omega_{\tilde{\cal A}_g}(d
\log \infty))$, cf.\ [F-C, p.195]. We denote the Chern classes of the Hodge
bundle by $\lambda_i$. The first Chern class $\lambda_1$ is ample on
${\cal A}_g$ and defines a morphism  $q:\tilde{\cal A}_g\to {\cal
A}_g^*$ to the Satake compactification. The $\Q$-subalgebra generated by the
$\lambda_i$ was determined in [vdG]. The proof there is given only for
$CH^*({\cal A}_g)$ and here we shall use no more than that.
\par
The Chow groups we discuss in this paper are always taken with
$\Q$-coefficients. We use equivariant classes on the level $\ell$-spaces,
cf.\ [E-G]. The classes of loci are always taken in the sense of
$Q$-classes, cf.\ [M3].  So if $Y\subset \tilde{\cal A}_g$ is a subvariety
then the class of $Y$ in $CH^*(\tilde{\cal A}_g)$ is taken with coefficient
$1/n$, where $n$ is the order of the automorphism group of the
completed (semi-)abelian variety corresponding to the generic point of $Y$.
The automorphisms must preserve the group structure of the 
semi-abelian part and extend to the completion. We shall frequently use the
following standard fact: if
$Y \subset X$ is a closed subvariety then we have an exact sequence
$$
CH_i(Y) \tto CH_i(X) \tto CH_i(X-Y) \tto 0.
$$
\par
We need several loci in the moduli space. We denote by $A_{g_1,\ldots,
g_r}$  the locus of products of principally polarized abelian varieties
of dimensions $g_1, \ldots, g_r$ with $\sum g_i=g$. There is a map
$$
{\cal A}_{g_1}\times \ldots \times{\cal A}_{g_r} \tto A_{g_1,\ldots,
g_r}  \subseteq {\cal A}_g
$$
of finite degree. The closure of $A_{g_1,\ldots,g_r}$ in a toroidal
compactification $\tilde{\cal A}_g$ is denoted by
$\tilde A_{g_1,\ldots,g_r}$.  We shall use the notation 
$$
\beta_t=q^{-1}(\coprod_{j\leq g-t} {\cal A}_j)
$$ 
for the locus of semi-abelian varieties with torus rank $\geq t$ on a
toroidal compactification  $\tilde{\cal A}_g$ and 
$$
\tilde{\cal A}_g^{(t)} = \tilde{\cal A}_g - \beta_{t+1}
$$
for the space of rank $\leq t$ degenerations. Sometimes we shall write
$\tilde{\cal A}_g^{(\geq t)}$ for $\beta_t$ and $\tilde{\cal
A}_g^{(\leq t)}$  for $\tilde{\cal A}_g^{(t)}$. 
\par
An important remark here is
that the space of rank
$1$-degenerations
$\tilde{\cal A}_g^{(1)}$ is canonical, i.e., does not depend on a choice of
polyhedral decomposition, cf.\ [M2]. The boundary $\tilde{\cal A}_g^{(1)}-
{\cal A}_g$ is a divisor
$D$. It admits a map of generic degree $2$ 
$$
j\colon {\cal X}_{g-1} \to D
$$
which identifies  $D$ with the universal family of Kummer varieties over
${\cal A}_{g-1}$. 
We have the following lemma, cf.\ [M2]:
\smallskip
\noindent
\proclaim (1.1) Lemma. The pull back  $j^*(D)$ of the $Q$-class of
$D$ restricted to a generic fibre $X$ equals $- \Theta$ with $\Theta$ the
theta divisor on $X$. \par
\smallskip
\noindent
\proclaim (1.2) Lemma. The cycle class $\lambda_1^{{g+1-t \choose 2}+1}$
vanishes on $ \beta_t$.
\par
\smallskip
\noindent
{\sl Proof.} This follows from the fact that $\lambda_1$ is ample on
${\cal A}_g^*$ and the dimension of ${\cal A}_{g-t}$. Indeed, a multiple
of $\lambda_1$ is represented by a hyperplane and the dimension of
$\beta_t$ is ${g+1-t \choose 2}$. \par
\vskip 2.0pc
The stack $\tilde{\cal A}_g$ is defined over $\Z$ and one can consider the
Chow ring and the subring generated by the Chern classes $\lambda_i$ of the
Hodge bundle $\E$  and the `boundary' classes  $\sigma_j$ for all fibres
$\tilde{\cal A}_g\otimes {\overline {\F}}_p$ and
$\tilde{\cal A}_g\otimes \C$. In this paper we restrict ourselves to the
complex case because we use references to papers where this restriction is
made. I see no serious obstacles for extending these results to all fibres.
\vskip 2.0pc
\noindent
{\bf \S 2 Ring Structure for} $g=1$ {\bf and} $g=2$.
\vskip1.5pc
\noindent
In this section we describe the Chow rings for $g=1$ and $g=2$. The
formulation that we give is intended to bring out similarities in ring
structure. The structure of the Chow ring of $\tilde{\cal A}_1$ is well
known. Recall that  $\lambda_1$ denotes the first Chern class of the Hodge
bundle; we denote the Q-class of the boundary
$\beta_1$ by $\sigma_1$. 
\par
\proclaim Theorem 1. The Chow ring of $\tilde{\cal A}_1$ is 
generated by $\lambda_1$ and $\sigma_1$  and is isomorphic to 
$$
\Q [ \lambda_1,\sigma_1]/(\lambda_1\sigma_1, \sigma_1-12\lambda_1)
\cong \Q [\lambda_1]/(\lambda_1^2).
$$
\par
Note that if we work over  $\Z$ the relation $\sigma_1\lambda_1$ should be
replaced by $2\sigma_1\lambda_1=0$ because $2\sigma_1$ represents the class
of a `physical' point. We thus expect over $\Z$ the ring $\Z[\lambda_1]/(24
\lambda_1^2)$, cf.\ the result of [E-G]. 
\vskip 1.0pc
The structure of the Chow ring of the moduli stack
$\overline{\cal M}_2$ of stable curves of genus $2$ was determined by Mumford
in [M3]. We can identify the stack $\overline{\cal M}_2$ with $\tilde{\cal
A}_2$ and obtain thus the structure of the Chow ring of $\tilde{\cal A}_2$.
Here
$\tilde{\cal A}_2$ is a canonical toroidal compactification and it was
already considered by Igusa. Because we use this result for $g=3$ we reprove
this result by a method we also use for $g=3$.
\par
For  every level $\ell$ there is a canonical compactification $\tilde{\cal
A}_2[\ell]$ of ${\cal A}_2[\ell]$, the moduli space of principally polarized
abelian surfaces with a full symplectic level $\ell$ structure. This
compactification is described in [Y]. The compactification is obtained by
adding a divisor $D[\ell]$ to ${\cal A}_2[\ell]$. The finite
group $G_{\ell}\cong Sp(2g, \Z / \ell\Z)$ acts on $\tilde{\cal A}_2[\ell]$
and there is an equivariant finite morphism
$p_{\ell}\colon\tilde{\cal A}_2[\ell] \tto  \tilde{\cal A}_2$. 
Classes in the Chow group $CH^k(\tilde{\cal A}_2[\ell])$ which are
$G_{\ell}$-invariant give rise to classes in $CH^k(\tilde{\cal A}_2)$. 
\par
We denote by $\sigma_i[\ell] \in CH^i(\tilde{\cal A}_2[\ell])$ the $i$-th
elementary symmetric function in the components of the divisor $D[\ell]$. 
Since it is $G_{\ell}$-invariant it defines a class $\sigma_i=\sigma_i[1]$
in $CH^i(\tilde{\cal A}_2)$ and this class has the property
$$
p_{\ell}^*(\sigma_i)= {\ell}^i\sigma_i[\ell]).
$$
Besides these elements we have the Chern classes $\lambda_i$ of the Hodge
bundle $\E$ on $\tilde{\cal A}_2$. In [vdG] I showed that these Chern
classes satisfy in the Chow ring $CH^*({\cal A}_g)$ the relation 
$$
(1+\lambda_1+\lambda_2)((1-\lambda_1+\lambda_2) =1
$$
on ${\cal A}_2$. The proof of the following theorem will show that this
relation still holds on $\tilde{\cal A}_2$.
\par
\proclaim
Theorem 2. The Chow ring of $\tilde{\cal A}_2$ is generated
by $\lambda_1,\lambda_2,\sigma_1$  and is isomorphic to 
$$
\Q [ \lambda_1, \lambda_2, \sigma_1]/I_2
$$
with $I_2$ the ideal generated by the relations
$$
\eqalign{
&(1+\lambda_1+\lambda_2)(1-\lambda_1+\lambda_2)=1, \cr
&\lambda_2\sigma_1=0,\cr
&\sigma_1^2= 22\sigma_1\lambda_1-120 \lambda_1^2.\cr
}
$$
The ranks of the Chow groups are: $1, 2, 2, 1$.\par
We can write this ring also as:
$$
\Q [ \lambda_1, \sigma_1]/(
\lambda_1^2\, \sigma_1 , (\sigma_1-10\, \lambda_1)(\sigma_1 - 12\, 
\lambda_1)). 
$$
The proof will consist of two steps: i) we bound the ranks of
$CH^k(\tilde{\cal A}_2)$ for $k=1,2,3$ and ii) by using intersection
numbers we show that the above generators are sufficient. The relations
then follow from the intersection numbers.
\par
We begin by bounding the ranks. There is a well-known finite morphism
$\A^3 \to {\cal A}_2-A_{1,1}$ which shows the triviality of the Chow groups
of ${\cal A}_2-A_{1,1}$. This shows that $CH^1(\tilde{\cal A}_2)$ is
generated by the class of $\tilde{A}_{1,1}$ and by $\sigma_1$. Since
$A_{1,1}$ has trivial Chow ring we see that $CH^2$ is generated by classes
living (in codimension $1$) on $D$. Now $D$ is the image of ${\cal
X}_1\to D$ and has in codimension $1$  two generators: the class of the zero
section and the class of a fibre. Since $\tilde{\cal A}_2$ is unirational we
find that $CH^3$ has rank $1$.
\par
\noindent
\proclaim (2.1) Conclusion. The ranks $r_k={\rm rank} ( CH^k(\tilde{\cal
A}_2))$ are bounded by $r_1\leq 2, r_2\leq 2$. Moreover $r_0=1, r_3=1$.\par
Now we calculate intersection numbers. We start with the intersection
numbers for the classes $\sigma_i$. One calculates
$$
\deg(\sigma_3[\ell])= \# {\hbox {\rm $0$-dim strata } } = 
{1\over 3} \ell \mu_1(\ell)\mu_2(\ell)= {1\over 12} \gamma/ \ell^3.
$$
leading to $\sigma_3= 1/12$ on the stack. Indeed, each of the $\mu_2(\ell)$
1-dimensional cusps carries $\mu_1(\ell)$ $0$-dimensional cusps and each of
these occurs on three transversally intersecting $1$-dimensional
cusps. A calculation (cf.\ [Y, T1])  shows that  for two different 
components
$F,G$ of $D[\ell]$ that intersect we have: $\deg(F^2G)=-2$. We thus find
$$
\deg(\sigma_1[\ell]\sigma_2[\ell]) = 3 \deg(\sigma_3([\ell])) -
2\ell\mu_1[\ell]\mu_2[\ell] = -{1\over 4}\gamma/\ell^3.
$$
Using proportionality (in cohomology) and the fact that $\lambda_1^2$ can be
represented by a compactly supported cycle on ${\cal A}_2$ (use the Satake
compactification) we find $\deg (\lambda_1^3)=(1/2880)\gamma$. Furthermore,
for 
$i: \tilde{\cal X}_1 \to {\tilde {\cal A}}_2$ we have
$$
i_*(1_{\tilde{\cal X}_1})= 2\sigma_1
$$
and   $i^*(\sigma_1) $ on $\tilde{\cal X}_1$ is in a fibre $F$ of ${\cal
X}_1 \to {\cal A}_1$ equal to $-2\Theta$ with $\Theta$ the theta divisor
there. On the other hand $i^*(\lambda_1)=[F]/24$,
where $F$ is a general fibre of $\tilde{\cal X}_1\to \tilde{\cal A}_1$. This
gives  $\deg (\lambda_1\sigma_1[\ell]^2)= (-1/24)\gamma/\ell^2$. The relation
$12\lambda_1=\sigma_1$ for $g=1$ gives 
$$
\sigma_2= 6\lambda_1\sigma_1
$$ 
for $g=2$.  We summarize these intersection numbers in a table.
\medskip
\vbox{
\noindent{\bf Table 2a $2\times 1$.}
\medskip\centerline{\def\quad{\hskip 0.6em\relax}
\def\quod{\hskip 0.5em\relax }
\vbox{\offinterlineskip
\hrule
\halign{&\vrule#&\strut\quod\hfil#\quad\cr
height2pt&\omit&&\omit&&\omit&\cr
&$2\backslash 1$&&$\lambda_1$&&$\sigma_1$&\cr
height2pt&\omit&&\omit&&\omit&\cr
\noalign{\hrule}
height2pt&\omit&&\omit&&\omit&\cr
&$\lambda_1^2$&&$1/2880$&&$0$&\cr
&$\lambda_1\sigma_1$&&$0$&&$-1/24$&\cr
height2pt&\omit&&\omit&&\omit&\cr
\cr } \hrule}
}}
\medskip
\par
\noindent
The table shows that $\lambda_1$ and $\sigma_1$ are linearly independent in
codimension $1$  as are $\lambda_1^2$ and $\lambda_1\sigma_1$ in codimension
$2$. Thus we have the generators. The restriction of $\lambda_2$ to any
boundary class is zero since the Hodge bundle possesses a line bundle
quotient on the closure $\tilde D$ of $D$.
\par
Moreover, $\deg(\lambda_1\lambda_2)=1/5760$ since $\lambda_1^2=2\lambda_2$
on ${\cal A}_2$ and this now shows that $\lambda_1^2=2\lambda_2$ on
$\tilde{\cal A}_2$. We also know that $\sigma_1^2$ is a linear combination of
$\lambda_1^2$ and $\lambda_1\sigma_1$ and can determine the coefficients if
we know
$\deg(\sigma_1^3)$. We know that $\sigma_1^2= -120 \lambda_1^2 + a\,
\lambda_1\sigma_1$ for some $a\in \Q$. In order to determine $a$ we shall
use the following lemma.
\par
\noindent
\proclaim (2.2) Lemma. The class of the products of elliptic curves
satisfies
$$
[\tilde{A}_{1,1}] = 5\lambda_1- {1\over 2}\sigma_1.
$$\par
This is well known, cf.\ [Ig]; compare also  [Mu2]. One can use the modular
form of weight $10$ which vanishes along $A_{1,1}$. Since
$12\lambda_1=\sigma_1$ in $g=1$ we get $(12\lambda_1-\sigma_1)(5\lambda_1
-(1/2)\sigma_1)=0$. This gives the relation 
$$
\sigma_1^2= 22\sigma_1\lambda_1 - 120 \lambda_1^2
$$
that we want and also the value for $\deg(\sigma_1^3)$.
\par
\medskip
\vbox{
\noindent{\bf  Table 2b. } The degree of the top $\sigma$'s.
\medskip\centerline{\def\quad{\hskip 0.6em\relax}
\def\quod{\hskip 0.5em\relax }
\vbox{\offinterlineskip
\hrule
\halign{&\vrule#&\strut\quod\hfil#\quad\cr
height2pt&\omit&&\omit&&\omit&\cr
&$\sigma_3$&&$\sigma_2\sigma_1$&&$\sigma_1^3$&\cr
height2pt&\omit&&\omit&&\omit&\cr
\noalign{\hrule}
height2pt&\omit&&\omit&&\omit&\cr
&$1/4$&&$-1/4$&&$-11/12$&\cr
height2pt&\omit&&\omit&&\omit&\cr
\cr } \hrule}
}}
\medskip
\noindent 
{\bf (2.3) Remark.} The $Q$-classes $\sigma_i$ with $i=1,\ldots$
satisfy $\sigma_1=[\beta_1]$, $\sigma_2=[\beta_{2}]$ here for $g=2$, but in
general $\sigma_t$ is not a multiple of $\beta_{t}$. 
\smallskip
Consider now the class $Z$ of the closure of the image cycle  of $j\colon
\tilde{\cal A}_1 \to \tilde{\cal A}_2$ obtained from sending $[X]$ to
$[X\times E]$, with
$E$ a fixed generic elliptic curve. The pull back of $\lambda_1$ is
$\lambda_1$ on $\tilde{\cal A}_1$ and the pull back of $\sigma_1$ is
$\sigma_1$. In other words, $12\, \lambda_1 -\sigma_1$ vanishes under pull
back. If we write
$Z= a \, \lambda_1^2 + b \, \lambda_1\sigma_1 $ this gives
$(12/2880)\, a - (-1/24) \, b =0$, i.e. $a=-10\, b$.
\smallskip
\noindent
\proclaim (2.4) Proposition. The $Q$-class of the closure of the
locus of trivial extensions of elliptic curves $[B_2]$  is
$$
[B_2]  =  120 \, \lambda_2 -\sigma_2 .
$$
\par
Note also that one has (cf.\ Mumford [M2] p. 368)
$$
[B_2] = N_0\, \sigma_1 = (5\lambda_1 - \sigma_1/2)\sigma_1.
$$
(cf.\ Mumford [M2] p. 368), and the two formulas agree by the relations
$\lambda_1^2=2\lambda_2$, $\sigma_1^2= 22\sigma_1\lambda_1 - 120
\lambda_1^2$ and $\sigma_2=6\lambda_1\sigma_1$.
\bigskip
\noindent
{\bf \S 3 Ring Structure for} $g=3$.
\bigskip
\noindent
We wish to describe the full Chow ring of the standard Delaunay-Voronoi
compactification for $g=3$. For a precise description of this
compactification see Tsushima [T1]. Tsushima constructs a smooth
compactification  $\tilde{\cal A}_3[\ell]$ for all $\ell\geq 3$. These
compactifications  are equipped with an action of $G_{\ell} = Sp(6, \Z /
\ell \Z)$ and with equivariant morphisms
$$
p_{\ell} \colon \tilde{\cal A}_3[\ell] \tto \tilde{\cal A}_3.
$$
Let  $ \gamma= \# G_{\ell}$ here be the order of this group. 
\par
The Chow ring $CH_{\Q}(\tilde{\cal A}_3)$  is defined to be the
invariant part of the Chow ring of $\tilde{\cal A}_3[\ell]$ in the sense
of [E-G].  It contains the $G_{\ell}$-invariant classes on $\tilde{\cal
A}_3[\ell]$. As a subring it contains the ring generated by the $\lambda_i$.
\par
The difference  $\tilde{\cal A}_3[\ell]-{\cal A}_3[\ell]$ is a union of
divisors. Let $\sigma_i[\ell] \in CH_{\Q}^i(\tilde{\cal A}_3[\ell])$ be the
$i$-th symmetric function of these divisors and set $\sigma_i=\sigma_i[1]$.
From the ramification properties along these divisors we derive
$$
p_{\ell}^*(\sigma_i) = {\ell}^i \sigma_i[\ell].
$$
\smallskip
\noindent
\proclaim Theorem 3. The Chow ring of $\tilde{\cal A}_3$ is generated by the
$\lambda_i, i=1,2,3$  and the $\sigma_i,  i=1,2$ and is isomorphic to the
graded ring (subscript is degree)
$$
\Q [ \lambda_1, \lambda_2,\lambda_3,  \sigma_1, \sigma_2]/I,
$$
with $I$ the ideal generated by the relations:
$$
\eqalign{
&(1+\lambda_1+\lambda_2 + \lambda_3)(1-\lambda_1 +\lambda_2
-\lambda_3)=1,\cr
&\lambda_3\sigma_1=\lambda_3\sigma_2=\lambda_1^2\sigma_2=0,\cr
&\sigma_1^3= 2016 \, \lambda_3 - 4 \, \lambda_1^2 \sigma_1 -24
\, \lambda_1\sigma_2 +{11\over 3}\, \sigma_2\sigma_1,\cr
&\sigma_2^2= 360\,  \lambda_1^3\sigma_1 -45 \, \lambda_1^2\sigma_1^2 
 +15 \, \lambda_1\sigma_2\sigma_1,\cr
&\sigma_1^2\sigma_2=3\, \sigma_2^2-30\, \lambda_1^2\sigma_1^2+2\, 
\lambda_1\sigma_1\sigma_2.\cr
}
$$
The ranks of the Chow groups are: $1, 2, 4, 6, 4, 2, 1$.
\par
For the proof we proceed as before. We first bound the ranks of the Chow
groups. Alternatively,  this could be done by using the results of Faber [F]
and the morphism $t:\overline{\cal M}_3 \tto \tilde{\cal A}_3$, see \S 4.
Analysis of the generators for $\overline{\cal M}_3$ leads to bounds on the
ranks of the Chow groups. But we proceed here more directly and use
parametrizations (as in [F]) to bound the ranks.
\par
\noindent
\proclaim (3.1) Proposition. The tautological ring of ${\cal A}_3$ generated
by the Chern classes $\lambda_i$ of the Hodge bundle is isomorphic to
$\Q[\lambda_1]/(\lambda_1^4)$.
\par
\noindent
{\sl Proof.} In [vdG] I showed that on ${\cal A}_3$ the  relations 
$$
\eqalign{
&(1+\lambda_1+\lambda_2+\lambda_3)(1-\lambda_1+\lambda_2-\lambda_3)=1,\cr
&\lambda_3=0\cr
}
$$ 
hold. This implies the result. $\square$
\smallskip
Recall that the locus of abelian varieties whose theta divisor is
singular yields a  divisor $N_0$ as defined by Mumford. Its class is given
by $N_0= 18\, \lambda_1-2\sigma_1$ on $\tilde{\cal A}_3$, cf.\ [I], [M2].
\par
Another divisor is obtained as the zero divisor of a modular form of weight
$140$. Let $\psi$ be the modular form of weight $140$ which is the product 
of the $8$-th powers of the 35
$\theta[\epsilon]$ for $\epsilon\neq 0$ and let $\Psi$ be its zero divisor
on $\tilde{\cal A}_3$. It is not difficult to see that the intersection of
$\Psi$ and $N_0$ is $240 A_{2,1}$. Analysis of how they intersect on $D$
leads to the following value for the class of $\tilde A_{2,1}$:
\par
\proclaim (3.2) Proposition. The $Q$-class of $\tilde A_{2,1}$ on $\tilde
{\cal A}_3$ is given by
$$ 
[\tilde A_{2,1} ] = {21\over 2}  \lambda_1^2 - {5\over 2} \lambda_1\sigma_1
+  {1\over 8} \sigma_1^2 + {1\over 24} \sigma_2.
$$
\par
\noindent
{\sl Proof.} Let $R$ be the $Q$-class of the closure of the locus of
semi-abelian varieties of $t$-rank $1$ whose abelian part is a product of
elliptic curves. We have $R= 5 \, \lambda_1 \sigma_1 - \sigma_2$ which comes
from genus $2$, cf.\ Lemma (2.2).  We have from [I] and [T1]:
$$ 
N_0\, \Psi = 240 \, [\tilde{A}_{2,1}] + 10\, R.
$$
We thus have
$$
(18\, \lambda_1 - 2\sigma_1)(140\lambda_1- 15 \, \sigma_1) = 240 \,
[\tilde A_{2,1}] + 50 \lambda_1\sigma_1 - 10 \sigma_{2}.
$$
\smallskip
\noindent
\proclaim (3.3) Proposition. The ring $CH^*(\tilde{\cal A}_3)$ is generated
by
$\lambda_1$ and the images of  classes on $\tilde D$.
\par
\noindent
{\sl Proof.} We know $[A_{2,1}]= (21/2)\lambda_1^2$ on ${\cal A}_3$. On the
other hand we have the Torelli map of degree $2$
$$
t: {\cal M}_3(\C)  \tto ({\cal A}_3 -A_{2,1})(\C) \qquad C \mapsto {\rm
Jac}(C)
$$
and from the fact that $CH^*({\cal M}_3)=
\Q[\lambda_1]/(\lambda_1^2)$ (cf.\ [F]) we see that $CH^*({\cal A}_3)$ is
generated by
$\lambda_1$ and $CH^*({\cal A}_{2,1})$. The description of the Chow
rings of ${\cal A}_1$ and ${\cal A}_2$ implies that $CH^*({\cal A}_{2,1})$
is generated by $p_2^*(\mu_1)$ (with $\mu_1$ the first Chern class of the
Hodge bundle on ${\cal A}_2$). But this class can be obtained from
$\lambda_1 [A_{2,1}]$. So $CH^*({\cal A}_3)$ is generated by
$\lambda_1$. $\square$
\smallskip
\noindent
\proclaim (3.4) Corollary. The ring $CH^*({\cal A}_3)$ is generated by
$\lambda_1$. The group $CH^1(\tilde{\cal A}_3)$  is generated by $\lambda_1$
and $\sigma_1$. These are linearly independent.
\par
\noindent
The linear independence follows e.g.\ from $\lambda_1^6\neq 0$ and
$\lambda_1^5\sigma_1=0$, which follows immediately by using the Satake
compactification.
\par
We now stratify the moduli space according to torus rank:
$$
\tilde{\cal A}_3 = {\cal A}_3 \cup \tilde{\cal A}_3^{t=1} \cup 
\tilde{\cal A}_3^{t=2} \cup \tilde{\cal A}_3^{t=3},
$$ 
with ${\rm codim}(\tilde{\cal A}_3^{t=j}) = j$.
\par
In view of  Proposition (3.3) we now first  determine the Chow groups of
$\tilde{\cal A}_3^{t\geq 1}=\tilde{\cal A}_3-{\cal A}_3$. From Section 1
we have a morphism $j: \tilde{\cal X}_2 \to \tilde{\cal
A}_3^{t\geq 1}$ with $\tilde{\cal X}_2$ the universal fibering over
$\tilde{\cal A}_2$ and this is of degree $2$ on ${\cal X}_2$. On $\tilde{\cal
X}_2$ we have the pull-backs $\tilde\lambda_i, \tilde\sigma_i$ of the classes
$\lambda_i, \sigma_i$  on $\tilde{\cal A}_2$ and the class $t$ coming from
the theta-divisor on the universal abelian surface. We normalize $t$ such
that $t|S$ is trivial, where $S$ is the zero-section. The class of $S$ is
denoted
$s$.
\smallskip
\noindent
\proclaim (3.5) Lemma. The Chow group $CH^1(\tilde{\cal X}_2)$ is generated
by
$t, \tilde\lambda_1$ and $\tilde\sigma_1$. 
\par
\noindent
{\sl Proof.} Consider the generic fibre $X_{\eta}$ with dual variety $\hat
X_{\eta}$. For the codimension $1$  Chow group we have $CH^1(X_{\eta}) =\Z
\oplus
\hat{X}_{\eta}$. But the only section of $\hat {\cal X}_2$ over ${\cal
A}_2$ is the zero section. If $L$ is a line bundle on $\tilde{\cal X}_2$
then $L \otimes t^m$ is trivial on $X_{\eta}$ for some $m$ and $\pi_*(L)$
is a line bundle on $\tilde{\cal A}_2$ with $\pi^*(\pi_*(L))\cong L$. So $L$
can be written as $L \cong t^{-m}
\otimes \pi^*(\lambda_1^a \otimes \sigma_1^b)$.
$\square$
\smallskip
Using the $g=2$ results and Lemma 1.1 we obtain the following table:
\medskip
\vbox{
\noindent{\bf  Table 3a.}  $4\times 1$ (i.e. codim $4$ $\times$ codim$1$)
on $\tilde{\cal X}_2$.
\medskip\centerline{\def\quad{\hskip 0.6em\relax}
\def\quod{\hskip 0.5em\relax }
\vbox{\offinterlineskip
\hrule
\halign{&\vrule#&\strut\quod\hfil#\quad\cr
height2pt&\omit&&\omit&&\omit&&\omit&\cr
&$4\backslash 1$&&$\tilde\lambda_1$&&$\tilde\sigma_1$&&$t$&\cr
height2pt&\omit&&\omit&&\omit&&\omit&\cr
\noalign{\hrule}
height2pt&\omit&&\omit&&\omit&&\omit&\cr
&$\tilde\lambda_1^3t$&&$0$&&$0$&&$1/1440$&\cr
&$\tilde\lambda_1^2 s$&&$1/2880$&&$0$&&$0$&\cr
&$\tilde\lambda_1\tilde\sigma_1 s$&&$0$&&$-1/24$&&$0$&\cr
height2pt&\omit&&\omit&&\omit&&\omit&\cr
\cr } \hrule}
}}
\medskip
\noindent
It follows from this that $CH^1(\tilde {\cal A}_3^{t\geq 1})$ is generated
by three elements: $\lambda_1\sigma_1, \sigma_1^2$ and $\sigma_2$. 
There is a relation between the generators of
$CH^2(\tilde{\cal A}_3)$ and of $CH^1(\tilde{\cal A}_3^{t\geq 1})$. By Lemma
1.1 we have:
$$
\sigma_1^2 \equiv -t  \bmod (\lambda_1\sigma_1)\qquad {\rm on}\quad
\tilde{\cal A}_3^{(t\leq 1)}.
$$
\proclaim (3.6) Proposition. The Chow group $CH^2(\tilde{\cal A}_3)$ is
generated by $4$ elements: $\lambda_1^2, \lambda_1\sigma_1,
\sigma_1^2$ and $\sigma_2$.
\par
Next we turn to $CH^3(\tilde{\cal A}_3)$. This group is generated by the
generator $\lambda_1^3$ of $CH^3({\cal A}_3)$ and the codimension $2$
classes on $\tilde{\cal A}_3^{t \geq 1}$. To determine these we decompose
the space $\tilde{\cal A}_3^{t\geq 1}$ in the following
four pieces (of dimensions $5, 4, 4, 3$):
$$
U:= {\tilde{\cal A}}_3^{t=1} - \pi_2^{-1}(A_{1,1})- S, \,\, \Xi:=
\pi_2^{-1}(A_{1,1}), \, 
 \tilde{\cal A}_3^{t\geq 2} \quad { \rm and }\quad S^0:=S-(S\cap
\pi^{-1}(A_{1,1})).  
$$
Note that $\Xi$ is a quotient of ${\rm Sym}^2({\cal X}_1)$ and it is fibred
over ${\rm Sym}^2({\cal A}_1)= \A^2$. It is a $\P^2$-bundle over $\A^2$,
from which we deduce $CH^1(\Xi)\cong \Q$, $CH^2(\Xi)\cong \Q$ and
$CH^k(\Xi)=0$ for $k\geq 3$. In accordance with Lemma (3.5) we have 
$CH^1(U)\cong \Q$.
\smallskip
\noindent
\proclaim (3.7) Lemma. The space $U$ has trivial Chow groups $CH^k$ for
$k\geq 2$. The Chow ring of $S^0$ is trivial (i.e. $\cong \Q$). Furthermore,
$CH^k(\Xi )\cong \Q$ for $k=1,2$  and $CH^k(\Xi )=0$  for $k>2$. 
The Chow ring of ${\cal A}_3^{t=1}$ is generated by $\tilde\lambda_1$ and
$t$.
\par
\smallskip
\noindent{\sl Proof.} As in [F] we can parametrize an open subset of
$\tilde{\cal A}_3^{t=1}$ by an open subset of the $\P^5-{\rm Quadric}$ of
plane quartics with a node.  Faber's Lemma (1.14) implies the first
statement. The second statement follows from the map $\A^3 \to S^0$, see
above. The  statements about $\Xi$ follow from the above. The last statement
follows from the decomposition $\tilde{\cal A}_3^{t=1} = U \cup \Xi \cup S$
and $S\cong {\cal A}_2$. $\square$
\smallskip
Summarizing, we find that $CH^3(\tilde{\cal A}_3)$ is generated by
$\lambda_1^3$, by the class of $S$, by the codimension $1$ class
$(\lambda_1t$) on $\Xi$ and by $CH^1(\tilde{\cal A}_3^{t\geq 2})$.      
\par
The morphism $\tilde{\cal X}_2 \tto \tilde{\cal A}_3^{t\geq 1}$
induces a morphism $\phi: \pi_2^{-1}(\tilde{\cal A}_2^{t\geq 1}) \tto
\tilde{\cal A}_3^{t\geq 2}$. The stratum $\tilde{\cal A}_3^{t=2}$ is now the
image of
$\pi_2^{-1}(\tilde{\cal A}_2^{t=1})$ which is fibred over $\tilde{\cal
A}_2^{t=1}$, which itself is a quotient of $\hat{\cal X}_1$. The fibres are
the image of a compactified $\G_m$-bundle over an elliptic curve. The stratum
$\tilde{\cal A}_3^{t=2}$ is thus the image of a compactified $\G_m$-bundle
over ${\cal X}_1\times_{{\cal A}_1}\hat{\cal X}_1$. We call the divisor
added to the $\G_m$-bundle $\Delta$. Since $\phi$ factors through ${\rm
Sym}^2({\cal X}/{\pm 1})$ we see that  our stratum ${\cal
A}_3^{t=2}-\phi(\Delta)$ is an affine bundle over a $\P^2$-bundle over
$\A^1$, the $j$-line, hence contributes only via its $CH^k\cong \Q$ for
$k=0,1,2$. The divisor
$\Delta$ is also a
$\P^2$-bundle over an affine line. Using these results on the  structure of
$\tilde{\cal A}_3^{t=2}$ we find:
\smallskip
\noindent
\proclaim (3.8) Lemma. For  $\tilde{\cal A}_3^{t=2}$ we have
${\rm rk}(CH^k) \leq 2$ for $k=1,2$ and  ${\rm rk}(CH^3) \leq 1$.
\par
\smallskip
\noindent
Collecting results, we get that $CH^2(\tilde{\cal A}_3^{t\geq 1})$ is
generated by the class of
$S$, the codimension $1$ class ($\tilde\lambda_1t$) on $\Xi$ and
$CH^1(\tilde{\cal A}_3^{t\geq 2})$. This gives as generators: $s$,
${\tilde\lambda_1}^2$,
$\tilde\lambda_1 \tilde\sigma_1$, $\tilde\lambda_1t$, and  $\tilde\sigma_1t$
and we have the bound $r_3 \leq 6$ for the rank $r_3$ of $CH^3(\tilde{\cal
A}_3)$.
\par
The $3$-dimensional toroidal stratum $\tilde{\cal A}_3^{t=3}$ can be viewed
using the map
$\pi_2$ and $\phi$ as a quotient under a finite map of a surface bundle over
a configuration of
$\P^1$'s. The fibres are a configuration of $\P^1\times \P^1$'s or a
configuration of $\P^2$'s and  $\P^2$'s blown up in three points, cf.\
Tsushima [T1, Lemma (4.4) and (7.1)]. Using this description one sees that
${\rm rk}(CH^1(\tilde{\cal A}_3^{t=3}))\leq2$. Indeed, we find three strata
that can contribute, but by writing down a rational function one gets a
non-trivial relation between them. For codimension $2$ we find
${\rm rk}(CH^2(\tilde{\cal A}_3^{t=3}))\leq 2$;  it is however not difficult
to see that the two 1-dimensional strata in $\tilde{\cal A}_3^{t=3}$
generate a $1$-dimensional space. But we do not need this latter estimate; it
suffices to remark that all possible generators of $CH^5(\tilde{\cal A}_3)$
lie in the ring generated by the $\lambda_i$ and the $\sigma_i$.
\par
Collecting again, we find that $CH^4(\tilde{\cal A}_3)$ is generated by a
generator of $CH^2(\Xi)$, by the two generators of $CH^2(\tilde{\cal
A}_3^{t=2})$ and by the generators of $CH^1(\tilde{\cal A}_3^{t=3})$.
This gives the bound  $r_4\leq 5$ (and similarly $r_5\leq 3$). Moreover, one
sees that the generators we gave can be expressed using the $\sigma$'s and
the $\lambda$'s.
\medskip
\vbox{
\noindent{\bf Table 3b } $3\times 2$ on $\tilde{\cal X}_2$. 
\medskip
\medskip\centerline{\def\quad{\hskip 0.6em\relax}
\def\quod{\hskip 0.5em\relax }
\vbox{\offinterlineskip
\hrule
\halign{&\vrule#&\strut\quod\hfil#\quad\cr
height2pt&\omit&&\omit&&\omit&&\omit&&\omit&&\omit&\cr
&$3\backslash 2$&&$\tilde\lambda_1^2$&&$\tilde\lambda_1\tilde\sigma_1$&&
$s$&&$\tilde\lambda_1t$&&$t\tilde\sigma_1$&\cr
height2pt&\omit&&\omit&&\omit&&\omit&&\omit&&\omit&\cr
\noalign{\hrule}
height2pt&\omit&&\omit&&\omit&&\omit&&\omit&&\omit&\cr
&$\tilde\lambda_1^3$&&$0$&&$0$&&$1/2880$&&$0$&&$0$&\cr
&$\tilde\lambda_1^2t$&&$0$&&$0$&&$0$&&$1/1440$&&$0$&\cr
&$\tilde\lambda_1s$&&$1/2880$&&$0$&&$1/5760$&&$0$&&$0$&\cr
&$\tilde\sigma_1s$&&$0$&&$-1/24$&&$0$&&$0$&&$0$&\cr
&$\tilde\lambda_1\tilde\sigma_1t$&&$0$&&$0$&&$0$&&$0$&&$-1/12$&\cr
height2pt&\omit&&\omit&&\omit&&\omit&&\omit&&\omit&\cr
\cr } \hrule}
}}
\smallskip
\noindent
\proclaim (3.9) Conclusion. The ranks $r_i={\rm rank}(CH^i(\tilde{\cal
A}_3^{t \geq 1}))$ are bounded by $r_0=r_5=1$, $r_1=3$, $r_2\leq 5$, $r_3\leq
5$ and $r_4\leq 3$. \par
\smallskip
\noindent
\proclaim (3.10) Conclusion. The ranks of the Chow groups of $\tilde{\cal
A}_3$ satisfy: $r_0=r_6=1$, $r_1=2$, $r_2\leq 4$, $r_3\leq 6$,
$r_4\leq 5$ and $r_5 \leq 3$. \par
\smallskip
\noindent
\proclaim (3.11) Lemma. We have $\lambda_3\sigma_i=0$ and
$\lambda_1^2\sigma_2=0$. \par
\noindent
{\sl Proof.} Since the Hodge bundle restricted (to a component of)
$\tilde D[\ell]$ possesses a subbundle of rank $2$ on each irreducible
component of $D[\ell]$ with trivial quotient its top Chern class vanishes
there. The second statement follows from Lemma (1.2).
$\square$. 
\bigskip
\smallskip
By the classical proportionality of Hirzebruch in cohomology we have  $\deg
\lambda_1^6 =
\gamma /181440$ in level $\ell$. We interpret this simply as 
$$
\lambda_1^6= {1\over 181440}
$$ 
in level $1$ omitting the degree. Note that $\tilde{\cal A}_3$ is known to
be unirational.
\par
The work of Tsushima implicitly contains all the information about the
intersection numbers on $\tilde{\cal A}_3$ and $\tilde{\cal D}$ that we
need. We give the results in the form of tables and indicate how they can
be obtained.
\par
We begin with intersection numbers for the sigma-classes, cf.\  Lemma
(7.7) of [T1] plus the correction in [T2].
\par
\medskip
\vbox{
\noindent{\bf Table 3c.} The $\sigma$'s. (Multiply entries by
$\gamma/\ell^6$ in level $\ell$.)
\medskip\centerline{\def\quad{\hskip 0.6em\relax}
\def\quod{\hskip 0.5em\relax }
\vbox{\offinterlineskip
\hrule
\halign{&\vrule#&\strut\quod\hfil#\quad\cr
height2pt&\omit&&\omit&&\omit&&\omit&&\omit&&
\omit&&\omit&&\omit&&\omit&&\omit&&\omit&\cr
&$\sigma_6$&&$\sigma_5\sigma_1$&&$\sigma_4\sigma_2$&&$\sigma_4\sigma_1^2$
&&$\sigma_3^2$&&$\sigma_3\sigma_2\sigma_1$&&$\sigma_3\sigma_1^3$
&&$\sigma_2^3$&&$\sigma_2^2\sigma_1^2$&&$\sigma_2\sigma_1^4$&&
$\sigma_1^6$&\cr
height2pt&\omit&&\omit&&\omit&&\omit&&\omit
&&\omit&&\omit&&\omit&&\omit&&\omit&&\omit&\cr
\noalign{\hrule}
height2pt&\omit&&\omit&&\omit&&\omit&&\omit
&&\omit&&\omit&&\omit&&\omit&&\omit&&\omit&\cr
&${1\over48}$&&${-1\over16}$&&${3\over16}$&&${13\over 48}$&&${41\over
144}$&&${1\over 16}$&&${-13\over 48}$&&${-15\over 16}$&&${-47\over
16}$&&${-445\over 48}$&&${ -4103\over 144}$&\cr
height2pt&\omit&&\omit&&\omit&&\omit&&\omit&&
\omit&&\omit&&\omit&&\omit&&\omit&&\omit&\cr
\cr} \hrule}}}
\bigskip
\medskip
\vbox{
\noindent{\bf  Table 3d. } Products $\lambda_1\times$ a product of 
$\sigma$'s.  (Multiply entries by $\gamma/\ell^5$ in level $\ell$.)
\medskip\centerline{\def\quad{\hskip 0.6em\relax}
\def\quod{\hskip 0.5em\relax }
\vbox{\offinterlineskip
\hrule
\halign{&\vrule#&\strut\quod\hfil#\quad\cr
height2pt&\omit&&\omit&&\omit&&\omit&&\omit
&&\omit&&\omit&&\omit&\cr
&$\lambda_1 \backslash
\sigma_{?}$&&$\sigma_5$&&$\sigma_4\sigma_1$&&$\sigma_3\sigma_2$&&$\sigma_3\sigma_1^2$
&&$\sigma_2^2\sigma_1$&&$\sigma_2\sigma_1^3$&&$\sigma_1^5$ &\cr
height2pt&\omit&&\omit&&\omit&&\omit&&\omit
&&\omit&&\omit&&\omit&\cr
\noalign{\hrule}
height2pt&\omit&&\omit&&\omit&&\omit&&\omit
&&\omit&&\omit&&\omit&\cr
&$\lambda_1$&&$0$&&$0$&&$1/48$&&$1/48$&&$-1/16$&&$-11/48$&&$-203/240$&\cr
height2pt&\omit&&\omit&&\omit&&\omit&&\omit
&&\omit&&\omit&&\omit&\cr
\cr } \hrule} }}
\medskip
\noindent
Another intersection number that we need is
$$
\lambda_1^3\sigma_1^3= {1 \over 720}.
$$
This follows easily from $\lambda_1^3|D= (1/2880)F$ with $F$ a generic fibre
of $\tilde D\to \tilde{\cal A}_2$ and the expression for $\sigma_1|D$ of
Lemma (1.1). Indeed, under the pull-back to $\tilde{\cal X}_2$ we have
$\deg(j^*(\lambda_1^3\sigma_1^2))=\deg(4t^2/2880)=1/360$. Dividing by the
degree of $j$ gives the result.

\smallskip
\noindent
The following tables can be constructed with this information. Together
with the bounds on the ranks of the Chow groups these will show that we found
the generators.
\medskip
\vbox{
\noindent{\bf Table 3e. } $1\times 5$.
\medskip\centerline{\def\quad{\hskip 0.6em\relax}
\def\quod{\hskip 0.5em\relax }
\vbox{\offinterlineskip
\hrule
\halign{&\vrule#&\strut\quod\hfil#\quad\cr
height2pt&\omit&&\omit&&\omit&\cr
&$1\backslash 5$&&$\lambda_1^5$&&$\lambda_1^3\sigma_1^2$&\cr
height2pt&\omit&&\omit&&\omit&\cr
\noalign{\hrule}
height2pt&\omit&&\omit&&\omit&\cr
&$\lambda_1$&&$1/181440$&&$0$&\cr
&$\sigma_1$&&$0$&&$1/720$&\cr
height2pt&\omit&&\omit&&\omit&\cr
\cr } \hrule}
}}
\medskip
\noindent
\bigskip
\noindent
\bigskip
\vbox{
\noindent{\bf Table 3f. } $2 \times 4$
\medskip
\medskip\centerline{\def\quad{\hskip 0.6em\relax}
\def\quod{\hskip 0.5em\relax }
\vbox{\offinterlineskip
\hrule
\halign{&\vrule#&\strut\quod\hfil#\quad\cr
height2pt&\omit&&\omit
&&\omit&&\omit&&\omit&\cr
&$2 \backslash
4$&&$\lambda_1^4$&&$\lambda_1^3\sigma_1$&&$\lambda_1^2\sigma_1^2$&&
$\lambda_1\sigma_2\sigma_1$&\cr
height2pt&\omit&&\omit&&\omit&&\omit&&\omit&\cr
\noalign{\hrule}
height2pt&\omit&&\omit
&&\omit&&\omit&&\omit&\cr
&$\lambda_1^2$&&$1/181440$&&$0$&&$0$&&$0$&\cr
&$\lambda_1\sigma_1$&&$0$&&$0$&&$1/720$&&$0$&\cr
&$\sigma_1^2$&&$0$&&$1/720$&&$0$&&$-11/48$&\cr
&$\sigma_2$&&$0$&&$0$&&$0$&&$-1/16$&\cr
height2pt&\omit&&\omit &&\omit&&\omit&&\omit&\cr
\cr } \hrule} }}
\medskip
\noindent
\medskip
\vbox{
\noindent{\bf  Table 3g.} $3\times 3$
\medskip
\medskip\centerline{\def\quad{\hskip 0.6em\relax}
\def\quod{\hskip 0.5em\relax }
\vbox{\offinterlineskip
\hrule
\halign{&\vrule#&\strut\quod\hfil#\quad\cr
height2pt&\omit&&\omit&&\omit&&\omit&&\omit&&\omit&&\omit&\cr
&$3\backslash 3$&&$\lambda_3$&&$\lambda_1^3$&&
$\lambda_1^2\sigma_1$&&$\lambda_1\sigma_2$&&$\lambda_1\sigma_1^2$&&
$\sigma_2\sigma_1$&\cr
height2pt&\omit&&\omit&&\omit&&\omit&&\omit&&\omit&&\omit&\cr
\noalign{\hrule}
height2pt&\omit&&\omit&&\omit&&\omit&&\omit&&\omit&&\omit&\cr
&$\lambda_3$&&$0$&&$1/1451520$&&$0$&&$0$&&$0$&&$0$&\cr
&$\lambda_1^3$&&$1/1451520$&&$1/181440$&&$0$&&$0$&&$0$&&$0$&\cr
&$\lambda_1^2\sigma_1$&&$0$&&$0$&&$0$&&$0$&&$1/720$&&$0$&\cr
&$\lambda_1\sigma_2$&&$0$&&$0$&&$0$&&$0$&&$0$&&$-1/16$&\cr
&$\lambda_1\sigma_1^2$&&$0$&&$0$&&$1/720$&&$0$&&$0$&&$-11/48$&\cr
&$\sigma_2\sigma_1$&&$0$&&$0$&&$0$&&$-1/16$&&$-11/48$&&$-47/16$&\cr
height2pt&\omit&&\omit&&\omit&&\omit&&\omit&&\omit&&\omit&\cr
\cr } \hrule}
}}
\medskip
{From} these tables and the bounds for the ranks we get the relations:
$$
\eqalign{ 
\sigma_3& = -40\, \lambda_1^2\sigma_1 +{44\over 3}
\lambda_1\sigma_2 -{1\over 3} \sigma_2\sigma_1, \cr
\sigma_1^3&= 2016 \, \lambda_3 - 4 \, \lambda_1^2 \sigma_1 -24
\, \lambda_1\sigma_2 +{11\over 3}\, \sigma_2\sigma_1.\cr
}
$$
Note that the class of the locus $\beta_3$ is known: $\lambda_1\sigma_2= {1
\over 4}\beta_{3}$. 
\bigskip
We wish to determine the class of $B_3$, that is, the $Q$-class of the
cycle which is the closure of the cycle $Z$ which is the image of the map
$$
j\colon {\cal A}_2 \tto {\cal A}_3, \qquad [X] \mapsto [X\times E],
$$
where $E$ is a generic elliptic curve. This is a map of degree $2$ in the
sense of stacks: $j_*(1)=2[Z]$. Under pull back we find
\par
\vbox{
\noindent
\medskip
\noindent{\bf Table 3h.} Intersection with $B_3$
\medskip\centerline{\def\quad{\hskip 0.6em\relax}
\def\quod{\hskip 0.5em\relax }
\vbox{\offinterlineskip
\hrule
\halign{&\vrule#&\strut\quod\hfil#\quad\cr
height2pt&\omit&&\omit&&\omit&&\omit
&&\omit&&\omit&&\omit&\cr
&&&$\lambda_3$&&$\lambda_1^3$&&$\lambda_1^2\sigma_1$&&
$\lambda_1\sigma_2$&&$\lambda_1\sigma_1^2$&&$\sigma_2\sigma_1$&\cr
height2pt&\omit&&\omit&&\omit&&\omit &&\omit&&\omit&&\omit&\cr
\noalign{\hrule}
height2pt&\omit&&\omit&&\omit&&\omit
&&\omit&&\omit&&\omit&\cr
&$\deg j^*$&&$0$&&$1/2880$&&$0$&&$0$
&&$-1/24$&&$-1/4$&\cr
height2pt&\omit&&\omit&&\omit&&\omit &&\omit&&\omit&&\omit&\cr
\cr } \hrule} }}
\medskip
\noindent
Note that because of $\lambda_2\sigma_1=0$ for $g=2$ we have in view of
Table 2.2  $j^*(\sigma_2\sigma_1) =
(3/11)j^*(\sigma_1^3)=(3/11)\times(-11/12)$. This gives:
\smallskip
\noindent
\proclaim (3.12) Proposition. The $Q$-class of the closure of the locus of
trivial extensions of abelian surfaces is given by
$$
[B_3]= 252\, \lambda_3 - 15\, \lambda_1^2\sigma_1 + 2 \lambda_1 \sigma_2.
$$
\par
We can also determine the class of $\tilde{A}_{1,1,1}$.
\smallskip
\noindent
\proclaim (3.13) Proposition. The $Q$-class of the closure of the locus of
products of three elliptic curves is given by
$$
[\tilde A_{1,1,1}]= -35\, \lambda_3 + {35\over 2}\, \lambda_1^3 -{25\over
4}\, \lambda_1^2\sigma_1  +{5\over 8} \, \lambda_1\sigma_2 + 
{5\over 8}\, \lambda_1\sigma_1^2 - {1\over 12} \,  \sigma_2\sigma_1
$$
\par
\noindent
{\sl Proof.} A priori we find (by restricting the Hodge bundle)
$$
[\tilde A_{1,1,1}] = 35(\lambda_1\lambda_2 - \lambda_3) + r
$$
with $r$ orthogonal to $\lambda_3$ and $\lambda_1^3$. We know that
$(12\lambda_1 -\sigma_1)$  annihilates the class of
$\tilde{A}_{1,1,1}$ and this suffices to determine $r$; alternatively  one
calculates the intersection numbers with our basis of $CH^3$ directly:
$$
[1/82944,1/13824,1/1152, 1/192, 1/96, 1/16]).
$$
\par
This implies  the assertions in Theorem 3 about the Chow groups $CH^i$ for
$i\leq 3$. As to $CH^4$,  we know that its rank is $\leq 5$. Moreover it
follows from the considerations above that $CH^4(\tilde{\cal A}_3)$ is
generated by elements of the subring generated by the $\lambda_i$ and the
$\sigma_i$. In view of the relations above this leaves us with the following
generators for $CH^4$:
$$
\lambda_1^4, \,  \lambda_1^3\sigma_1, \, \lambda_1^2\sigma_1^2, \,
\lambda_1\sigma_1\sigma_2, \, \sigma_1^2\sigma_2, \, \sigma_2^2.
$$
I claim that we have two independent relations between these. 
Indeed, all these classes are represented by cycles on 
$\tilde{\cal A}_3^{t\geq 2}$ and as we
observed above we have $CH^2(\tilde{\cal A}_3^{t\geq 2})\leq 4$. Using the
intersection numbers we find the relations expressing $\sigma_2^2$ and
$\sigma_1^2\sigma_2$ in terms
of the others. Similarly, in codimension $5$ the classes generating the Chow
group lie in the subring generated by the $\lambda_i$ and the $\sigma$'s.
But the codimension $5$ part of this subring has rank $2$. So we get
$r_5=2$. The unirationality of
$\tilde{\cal A}_3$ implies $r_6=1$. The completes the proof of Theorem~3.
\vskip 1.0pc
\noindent
{\bf \S 4 Comparison with the Chow ring of $\overline{\cal M}_3$}
\bigskip
\noindent
Let $\overline{\cal M}_3$ be the moduli space of stable curves of genus $3$.
The usual Torelli map $t:{\cal M}_3 \to {\cal A}_3$ can be extended to a
morphism $\bar t:\overline{\cal M}_3 \tto \tilde{\cal A}_3$ in the following
way: associate to a curve $C$ the (coarse) moduli space $P$ of semi-stable
torsion free rank $1$ sheaves of degree $g-1$ on $C$. This moduli space comes
provided with a `theta' divisor $B$ whose points correspond to the
semi-stable sheaves with $h^0 \neq 0$. Such a pair $(P,B)$ defines in a
natural way a point of the Delaunay-Voronoi compactification $\tilde{\cal
A}_3$. The Delaunay decomposition of $\Z^t$ with $t=$ torus rank of $P$
corresponding to $C$ is obtained as follows. Let $\Gamma$ be the dual
graph of the curve. The space $C_1(\Gamma)$ is provided with the standard
Euclidean Delaunay decomposition, i.e. the standard cube, its faces and
translates. It induces a decomposition on the linear subspace
$H_1(\Gamma)$.  This is the Delaunay decomposition of $\Z^t$, cf.\ [A, AN].
\par
The Chow group $CH^1(\overline{\cal M}_3)$ is generated by $t^*(\lambda_1)$,
$t^*(\sigma_1)=\delta_0$, and a new element denoted in [F] by
$\delta_1$. It comes from the reducible curves with components of genus $1$
and  genus $2$.  The locus $\Delta_1$ of these reducible curves
maps to the codimension $2$ locus $\tilde A_{2,1}$. So in 
codimension $1$ the map
$t_*: CH^1(\overline{\cal M}_3) \tto CH^1(\tilde{\cal A}_3)$ is given by
$$
\eqalign{
t_*(\lambda_1)&= 2\lambda_1,\cr
t_*(\delta_0)&= 2\sigma_1,\cr
t_*(\delta_1)&=0.\cr
}
$$
In codimension $2$ the group $CH^2(\overline{\cal M}_3)$ is generated by the
four pull backs under $t^*$ of the classes $\lambda_1^2, \lambda_1\sigma_1,
\sigma_1^2$ and $\sigma_2$ and three new generators. These are
$$
\lambda_1\delta_1, 
\sigma_1\delta_1, \delta_{1,1},
$$
where $\delta_{1,1}$ is the $Q$-class of the closure of the locus of genus
$3$ curves which are a string of elliptic curves, i.e. with dual graph
\def\lijn{\hskip-1mm\raise0.9mm\hbox{\vrule height0.2pt
width15mm}\hskip-1mm}
$${\buildrel 1 \over \bullet}\lijn {\buildrel 1 
\over \bullet}\lijn{\buildrel 1 \over \bullet}
$$
\noindent
These elements are mapped to zero under $t_*$ and with the notation $[a,b,c,d]$ for 
$a\lambda_1^2+b\lambda_1\sigma_1+c\sigma_1^2+d\sigma_2$ the images under 
the map $t\colon {\overline{\cal M}}_3 \tto \tilde{\cal A}_3$ are:
$$
\eqalign{
t_*(\delta_{00}) & = 2\times[0,0,0,1] \cr
t_*(\xi_0) & = 2\times[0,4,-1,1] \cr
t_*(\xi_1)&= 2\times[0,5/2,0,-1/2] \cr
t_*(\eta_1)& = 2\times [63/2,-15/2,3/8,1/8] = 6\times [\tilde A_{2,1}]\cr
}
$$
For the reader's convenience we also give
$$
\eqalign{
t_*(\kappa_2)&= 2\times[41/2,-7/2,1/8,1/24]\cr
t_*(\delta_1^2)&= 2\times[-21/2,5/2,-1/8,-1/24]= -2\times [\tilde
A_{2,1}]. \cr }
$$
We have (cf.\ [M2]):
$t_*(\xi_0+2\xi_1)=2\times (9\, \lambda_1-\sigma_1)\, \sigma_1 =[N_0]$.
{From} [F, p.\ 368] we have the relation
$$
\delta_1^2= 3 \lambda_1\delta_1 - {1 \over 3}\sigma_1\delta_1 - {21\over 2}
\lambda_1^2 +{5\over 2}\lambda_1\sigma_1 - {1 \over 8} \sigma_1^2 - {1
\over 24} \sigma_2.
$$
Under $t_*:CH^3(\overline{\cal M}_3)\tto \tilde{\cal A}_3$ the images of
the classes $[(b)]_Q, [(d)]_Q, [(f)]_Q, [(g)]_Q$ are zero, while the image
of $[(i)]_Q$ is $[\tilde A_{1,1,1}]$.
\par
The following table gives the images of the generators of
$CH^3(\overline{\cal M}_3)$ under $t_*$.
\par
\medskip
\vbox{
\noindent{\bf  Table 4a.} Image under $t_*$
\medskip
\medskip\centerline{\def\quad{\hskip 0.6em\relax}
\def\quod{\hskip 0.5em\relax }
\vbox{\offinterlineskip
\hrule
\halign{&\vrule#&\strut\quod\hfil#\quad\cr
height2pt&\omit&&\omit&&\omit&&\omit&&\omit&&\omit&&\omit&\cr
&$t_*(\cdot )/2 $&&$\lambda_3$&&$\lambda_1^3$&&
$\lambda_1^2\sigma_1$&&$\lambda_1\sigma_2$&&$\lambda_1\sigma_1^2$&&
$\sigma_2\sigma_1$&\cr
height2pt&\omit&&\omit&&\omit&&\omit&&\omit&&\omit&&\omit&\cr
\noalign{\hrule}
height2pt&\omit&&\omit&&\omit&&\omit&&\omit&&\omit&&\omit&\cr
&$t_*(\lambda_3)/2$&&$1$&&$0$&&$0$&&$0$&&$0$&&$0$&\cr
&$t_*(\lambda_1^3)/2$&&$0$&&$1$&&$0$&&$0$&&$0$&&$0$&\cr
&$ t_*[(a)]_Q/2$&&$0$&&$0$&&$0$&&$4$&&$0$&&$0$&\cr
&$ t_*[(c)]_Q/2$&&$0$&&$0$&&$30$&&$-6$&&$0$&&$0$&\cr
&$ t_*[(e)]_Q/2$&&$0$&&$0$&&$-40$&&$32/3$&&$0$&&$-1/3$&\cr
&$t_*[(h)]_Q/2$&&$0$&&$0$&&$-25/2$&&$-5/2$&&$-5/4$&&$1/4$&\cr
&$t_*[(i)]_Q/2$&&$-35/2$&&$35/4$&&$-25/8$&&$5/16$&&$5/16$&&$-1/24$&\cr
&$t_*(\eta_0)/2$&&$756$&&$0$&&$-45$&&$6$&&$0$&&$0$&\cr
height2pt&\omit&&\omit&&\omit&&\omit&&\omit&&\omit&&\omit&\cr
\cr } \hrule}
}}
\medskip
A check on this is given by the relation
$$
504\, \lambda_3 = {1\over 2} \delta_{000} + \xi_{01} + {2 \over 3} \eta_0
$$
which can be deduced from [F], cf.\ [F2], p.\ 77. Here $\xi_{01}$ is the
locus of curves  consisting of an elliptic curve which intersects a
degenerate elliptic curve in $2$ points. Under $t_*$ this is translated into
$ 504 \,
\lambda_3 = 2\lambda_1\sigma_2 + t_*([(c)]_Q) + 2B_3$. So by the expression
for $[B_3]$ we find $ t_*(\xi_{01})= 6\lambda_1
(5\lambda_1\sigma_1-\sigma_2) $ which fits with the formula for $[(c)]_Q$
and the one for $[A_{1,1}]$  in genus $2$. Another check is obtained from
the formula
$$
\delta_0^3 = {47\over 15}[(a)]_Q-{54\over 5}[(b)]_Q-{54\over 5}[(c)]_Q
-{89\over 15}[(d)]_Q -11[(e)]_Q+{8\over 5}[(f)]_Q+8[(g)]_Q +{8 \over
3}\eta_0
$$
(cf [F,p.\ 411]). Under $t_*$ one finds our relation for $\sigma_1^3$.
\bigskip
\noindent
{\bf \S 5 Concluding Remarks.}
\bigskip
\noindent
In general one cannot expect that the Chow ring of a suitable
compactification $\tilde{\cal A}_g$ is generated by the tautological
classes $\lambda_i$ and $\sigma_i$. But it makes sense to consider the
subring generated by these elements. The relations satisfied by the
$\lambda_i$ are known. I do not know whether  a relation of the sort:
$$
\sigma_1^g = \zeta(1-2g) \lambda_g + \quad {\hbox{{\rm classes on }
$\tilde{\cal A }_g^{t\geq 2}$ }}
$$
holds in the Chow ring. In cohomology it does. 
Another approach would be to study the Chow rings of the canonical partial 
compactification $\tilde{\cal A}_g^{(1)}$. For $g=1, 2 $ and $3$ we find 
the following rings:
\par
\noindent
$g=1$:
$$
\Q [ \lambda_1,\sigma_1]/(\lambda_1\sigma_1, \sigma_1-12\lambda_1)
$$
$g=2$:
$$
\Q[\lambda_1,\sigma_1]/(\lambda_1\sigma_1,\sigma_1^2+120\, \lambda_1^2).
$$
$g=3$:
$$
\Q[\lambda_1,\lambda_3,
\sigma_1]/(\lambda_1^4-8\lambda_3\lambda_1,\lambda_1^2\sigma_1,\sigma_1^3-2016\,
\lambda_3,\lambda_3\sigma_1).
$$
\par
As a final remark we observe that the Chow ring of $\tilde{\cal X}_2$ is
generated as an algebra over  the Chow ring  of $\tilde{\cal A}_2$ (using the
map
$\pi_2^*$) by the  element $t$. It satisfies the quadratic equation
$$
t^2= -\tilde\lambda_1^2 + 2\, s + t\tilde\sigma_1/8.
$$
Indeed, a priori, $t^2$ is a linear combination 
$a\tilde\lambda_1^2+b\tilde\lambda_1 \tilde\sigma_1+cs+
d \tilde\lambda_1 t+ et \tilde\sigma_1$. But in view of the $3\times
2$-table  for $\tilde{\cal X}_2$ the relation  $ts=0$ implies that
$b=0$. Because  $s^2=\tilde\lambda_2s$ and $t^2$ has degree 
$2$ in a general fibre we find $a=-1$ and $c=2$. From [vdG] we know 
$(\pi_2)_*(t^3/3!)=\sigma_1/24$. This implies $d=0$ and $e=1/8$.
\par
\bigskip
\bigskip
\noindent
{\bf References.}
\bigskip
\noindent
[A] V.\ Alexeev: Compactified Jacobians. Preprint, {\tt alg-geom 9608012}
\smallskip
\noindent
[AN] V.\ Alexeev, I.\ Nakamura: On Mumford's construction of degenerating
abelian varieties. Preprint, {\tt alg-geom 9608014}
\smallskip
\noindent
[E-G] D.\ Edidin, W.\ Graham: Equivariant intersection theory. Preprint,
{\tt alg-geom 9603008} 
\smallskip
\noindent
[F] C. Faber: Chow rings of moduli spaces of
curves. I: The Chow ring of
$\overline{\cal M}_3$. {\sl Ann.\ of Math.\ \bf 132  }, p. 331-419 (1990)
\smallskip
\noindent
[F2] C. Faber: Intersection-theoretical computations on $\overline{\cal
M}_g$. In: {\sl Parameter Spaces}, Banach Center Publ.\ {\bf 36}, Warszawa
1996.
\smallskip
\noindent
[F-C] G.\ Faltings, C-L.\ Chai: Degeneration of abelian varieties.
Ergebnisse der Math. 22. Springer Verlag 1990.
\smallskip
\noindent
[vdG] G. van der Geer: Cycles on the moduli space of abelian varieties.
Preprint. {\tt alg-geom 9605011}
\smallskip
\noindent
[H] F.\ Hirzebruch: Automorphe Formen und der Satz von Riemann-Roch.
In: {\sl Sympo\-sium Internacional de Topologia Algebrica} (M\'exico 1956),
129-144. La Universidad Aut\'o\-no\-ma de M\'exico 1958 (= Collected Papers
I, Springer Verlag, p.\ 345.)
\smallskip
\noindent
[I] J.\ Igusa: A desingularization problem in the theory of Siegel modular
functions. {\sl Math.\ Annalen \bf 168}, p.\ 228-260 (1967).
\smallskip
\noindent
[M1] D.\ Mumford: Hirzebruch's Proportionality Theorem in the
non-compact case. {\sl Inv.\ Math.\ \bf 42} (1977) , 239-272.
\smallskip
\noindent
[M2] D.\ Mumford: On the Kodaira dimension of the Siegel modular variety.
In {\sl Open Problems, \rm Eds C.\ Ciliberto, F.\ Ghione, F.\
Orecchia}, Springer Lecture Notes 997 , p. 348-375, Springer Verlag
1983.  
\smallskip
\noindent
[M3] D.\ Mumford: Towards an enumerative geometry of the moduli
space of curves. In: {\sl Arithmetic and Geometry} (dedicated to
I.\ Shafarevich) Vol II. Progress in Mathematics 36, Birkh\"auser
1983.
\smallskip
\noindent
[N] Y.\ Namikawa: A new compactification of the Siegel space and
degenerations of abelian varieties I,II. {\sl Math.\ Ann.\ \bf 221} (1976,
p.\ 97-141, 201-241.
\smallskip
\noindent
[T1] R. Tsushima: A formula for the dimension of spaces of Siegel
modular cusp forms. {\sl Amer.\ Journal of Math.\ \bf 102} (1980), p.\
937-977 
\smallskip
\noindent
[T2] R. Tsushima: On the spaces of spaces of Siegel
modular cusp forms. {\sl Amer.\ Journal of Math.\  \bf 104} (1982), p.\
843-885 
\smallskip
\noindent
[Y] T. Yamazaki: On Siegel modular forms of degree two. {\sl American
Journal of Math.\ \bf 98}, p. 39-53 (1976)
\bigskip
\par
Gerard van der Geer
\par
Faculteit WINS
\par
Universiteit van Amsterdam
\par
Plantage Muidergracht 24
\par
1018 TV Amsterdam
\par
Nederland
\par
e-mail: {\tt geer@wins.uva.nl}
\end